\documentclass[letterpaper, 10 pt, conference]{ieeeconf}  % Comment this line out
\IEEEoverridecommandlockouts                              % This command is only
\usepackage{graphicx} % for pdf, bitmapped graphics files
\usepackage{amsmath} % assumes amsmath package installed
\usepackage{amssymb}  % assumes amsmath package installed
\usepackage{amstext,enumerate}
\usepackage{subfigure}
\usepackage{url}
\usepackage[symbol]{footmisc}
\usepackage{algorithm}
\usepackage{algpseudocode}

\usepackage{latexsym, color,booktabs}

\newtheorem{assum}{Assumption}
\newtheorem{prop}{Proposition}
\newtheorem{thm}{Theorem}
\newtheorem{rem}{Remark}
\newtheorem{lem}{Lemma}
\newtheorem{cor}{Corollary}
\newtheorem{defn}{Definition}

\title{\LARGE \bf
Entropy for incremental stability of nonlinear systems \\ under disturbances  
}
\author{Michelle S. Chong
\thanks{M. Chong is with the Control Systems Technology section at the Department of Mechanical Engineering, Eindhoven University of Technology. 
        {\tt\small m.s.t.chong@tue.nl} }
}

\begin{document}
\maketitle
\thispagestyle{empty}
\pagestyle{empty}

%%%%%%%%%%%%%%%%%%%%%%%%%%%%%%%%%%%%%%%%%%%%%%%%%%%%%%%%%%%%%%%%%%%%%%%%%%%%%%%%
\begin{abstract}
Entropy notions for $\varepsilon$-incremental practical stability and incremental stability of deterministic nonlinear systems under disturbances are introduced. The entropy notions are constructed via a set of points in state space which induces the desired stability properties, called an approximating set. We provide conditions on the system which ensures that the approximating set is finite. Lower and upper bounds for the two estimation entropies are computed. The construction of the finite approximating sets induces a robust state estimation algorithm for systems under disturbances using quantized and time-samples measurements. 
\end{abstract}

%%%%%%%%%%%%%%%%%%%%%%%%%%%%%%%%%%%%%%%%%%%%%%%%%%%%%%%%%%%%%%%%%%%%%%%%%%%%%%%%
\section{Introduction} \label{sec:intro}

Entropy quantifies the rate at which a dynamical system generates information. Its role in feedback control 
translates to the amount of information needed by the controller in order to affect a desired behaviour of the plant. This context is most applicable when the measurements of the plant is only available intermittently such as in the case of network control systems where the plant and controller communicate over a channel with finite capacity, which is a problem in consideration for over 30 years, see \cite{nair2007feedback} for a survey. 

The notion of entropy provides an abstract approach to characterise the amount of data (control law) needed to achieve a desired behaviour (stabilization or invariance) without being concerned about \textit{how} the data (control law) is generated. The pioneering work by Nair et. al. \cite{nair2004topological} introduced the notion of topological feedback entropy for discrete-time systems, which counts the number of open covers in state space. Another approach was introduced in \cite{tatikonda2004control} in discrete-time and \cite{colonius2009invariance} in continuous-time which counts the `spanning set' of open-loop control functions to achieve invariance. This was followed by \cite{colonius2012minimal} and \cite{colonius2021entropy} for exponential and practical stabilization, respectively.

We take the same approach in characterising entropy for incremental stability, which compares arbitrary solutions of a dynamical systems with themselves, instead of an equilibrium point or a particular trajectory \cite{zamani2011backstepping}. The incremental stability property has been identified to be applicable to state estimation, synchronisation of dynamical systems and in finite abstractions for nonlinear systems, to name a few. See \cite{zamani2011backstepping} for a historical account and a list of applications. Here, we are motivated by the state estimation problem using quantized and time-sampled measurements. 

The usage of entropy for state estimation under finite data rates has been studied in \cite{savkin2006analysis, nair2013nonstochastic, matveev2016observation} for discrete-time systems and in \cite{liberzon2017entropy,matveev2019observation, kawan2021remote} for continuous-time systems. With the exception of \cite{nair2013nonstochastic} which considered \textit{linear} systems under disturbances, we are not aware of any other works that consider \textit{nonlinear} systems in the presence of disturbances. Under this setup, standard state estimation algorithms where the measurement is continuously available typically do not yield incrementally stable properties, as the estimate typically converges to a neighbourhood of the true state. This leads to the notion of $\varepsilon$-incremental \textit{practical} stability considered in this paper, which has ties to the incremental input-output-to-state stability (i-IOSS) property and its application to state estimation \cite{sontag1997output}, but is not pursued here for conciseness.

In this paper, we introduce and compute bounds for new estimation entropy notions for $\varepsilon$-incremental practical stability and incremental stability. This extends the estimation entropy notion in \cite{liberzon2017entropy}, where only incremental exponential stability is considered. Like in \cite{liberzon2017entropy}, we require estimates to converge at a prescribed rate characterised by a class $\mathcal{KL}$ function (which includes the exponential rate considered in \cite{liberzon2017entropy}) and to a desired accuracy $\varepsilon>0$. We approach the estimation entropy notion by constructing sets of points of the state space, which lead to the desired incremental stability properties. This approach follows the classical construction of entropy for dynamical systems by Bowen \cite{bowen1971entropy} and Dinaburg \cite{dinaburg1971relations}, which was used in \cite{colonius2009invariance, colonius2012minimal,colonius2021entropy} for stabilization and in \cite{savkin2006analysis, liberzon2017entropy} for estimation, where the sets are known as `spanning sets', which we call `approximating sets' here. To ensure the existence of finite approximating sets, we needed to modify the incremental stability notions that can be achieved with the approximating sets, see Remark \ref{rem:def_exist}. We obtain the upper bounds on our notions of estimation entropy using matrix measures (also used in \cite{liberzon2017entropy}), and the lower bound on each estimation entropy notions is obtained via a volume growth argument.

The construction of the finite approximating sets induces a state estimation algorithm for systems under disturbances using quantized and time-sampled measurements. The proposed algorithm is a modification of the iterative procedure in \cite{liberzon2017entropy} to account for estimation inaccuracy due to the presence of disturbances. We provide convergence guarantees in that the estimate converges to a neighbourhood of the true state where the size of the neighbourhood depends on the  estimation inaccuracy $\varepsilon$.

The paper is organised as follows. We start with the preliminaries in Section \ref{sec:prelim} and formulate the problem in Section \ref{sec:problem}. The notions of estimation entropy for $\varepsilon$-incremental practical stability and incremental stability are defined in Section \ref{sec:entropy} via approximating sets. System properties which lead to the existence of \textit{finite} approximating sets are provided in Section \ref{sec:finite_sets}. We compute bounds on each notions of estimation entropy in Section \ref{sec:bounds}. We then propose a robust state estimation scheme using quantized and time-sampled measurements in Section \ref{sec:algo} and conclude the paper with Section \ref{sec:conclude}. Proofs of all results are provided in the Appendix.

\section{Preliminaries} \label{sec:prelim}
 \subsection{Notation}
 \begin{itemize}
 	\item Let $\mathbb{R}=(-\infty,\infty)$, $\mathbb{R}_{\geq 0}=[0,\infty)$, $\mathbb{R}_{>0}=(0,\infty)$. 
 	\item Let $\mathbb{N}_{\geq i}=\{i,i+1,i+2,\dots\}$. A finite set of integers $\{i,i+1,i+2,\dots,i+k\}$ is denoted as $\mathbb{N}_{[i,i+k]}$.
%	\item The number of $k$-element subsets of an $n$-element set is denoted $\binom{n}{k}$. 
 	%\item Let $(u,v)$ where $u\in\mathbb{R}^{n_u}$ and $v\in\mathbb{R}^{n_v}$ denote the column vector $(u^T,v^T)^{T}$.
% 	\item The cardinality of a set $\mathcal{J}$ is denoted as $\#\left(\mathcal{J}\right)$.
% 	\item The identity matrix of dimension $n$ is denoted by $\mathbb{I}_{n}$. % and a matrix of dimension $m$ by $n$ with all elements $1$ is denoted by $\mathbf{{1}}_{m\times n}$.
%	\item A diagonal matrix with matrices $d_i$, $i\in\mathbb{N}_{[1,n]}$ is denoted by $\textrm{diag}(d_1,d_2,\dots,d_n)$.
% 	\item The Euclidean norm of a vector $x \in \mathbb{R}^{n}$, is denoted $|x|:=\sqrt{x^{T}x}$ and for a matrix $A\in\mathbb{R}^{n\times n}$, $|A|:=\sqrt{\lambda_{\max}(A^T A)}$.
    \item The determinant of a matrix $A$ is denoted by $\textrm{det}(A)$ and the sum of its diagonal entries as $\textrm{tr}(A)$.
% 	\item Given a symmetric matrix $P$, its maximum (minimum) eigenvalue is denoted by $\lambda_{\max}(P)$ $(\lambda_{\min}(P))$.
 	\item For a given vector $x\in\mathbb{R}^{n}$ and matrix $A\in\mathbb{R}^{n\times n}$, let $|x|$ and $\|A\|$ denote a chosen norm in $\mathbb{R}^{n}$ and $\mathbb{R}^{n\times n}$, respectively.  In this paper, we find it convenient to use the infinity norm $|x|:= \underset{i\in\mathbb{N}_{[1,n]}}{\max} \left| x_i \right|$ and the induced matrix norm on $\mathbb{R}^{n\times n}$ corresponding to the chosen norm on $\mathbb{R}^{n}$ is $|A|:= \underset{i\in\mathbb{N}_{[1,n]}}{\max} \underset{j\in\mathbb{N}_{[1,n]}}{\sum}|a_{ij}|$, where $a_{ij}$ is the row $i$-th and column $j$-th element of matrix $A$. 
	\item Given a point $x\in\mathbb{R}^{n}$, the closed ball with radius $\delta$ around $x$ is denoted as $B(x,\delta)=\{ z \in \mathbb{R}^{n} \, |\, |z-x|\leq \delta\}$. When the infinity norm is used, the closed ball $B(x,\delta)$ is a hypercube.
	
%	\item Given a non-empty closed set $\mathcal{A}\subset \mathbb{R}^n$, $|x|_{\mathcal{A}}:=\underset{y\in\mathcal{A}}{\inf}|x-y|$.
     
%    \item Let $\mathcal{B}^{\mathcal{A}}_{\Delta}:=\left\{ x\in\mathbb{R}^{n} : |x|_{\mathcal{A}} \leq \Delta \right\}$.
 	\item A continuous function $\alpha:\mathbb{R}_{\geq 0}\to\mathbb{R}_{\geq 0}$ is a class $\mathcal{K}$ function, if it is strictly increasing and $\alpha(0)=0$; additionally, if $\alpha(r)\to\infty$ as $r\to\infty$, then $\alpha$ is a class $\mathcal{K}_{\infty}$ function. A continuous function $\beta:\mathbb{R}_{\geq0}\times \mathbb{R}_{\geq 0} \to \mathbb{R}_{\geq 0}$ is a class $\mathcal{KL}$ function, if: (i) $\beta(.,s)$ is a class $\mathcal{K}$ function for each $s\geq 0$; (ii) $\beta(r,.)$ is non-increasing and (iii) $\beta(r,s)\to 0$ as $s\to \infty$ for each $r\geq 0$.
\end{itemize}

\subsection{Matrix measures}
For a given matrix $A\in\mathbb{R}^{n\times n}$, the matrix measure $\mu:\mathbb{R}^{n \times n} \to \mathbb{R}$ is the one-sided derivative of the induced norm at $I\in\mathbb{R}^{n\times n}$ in the direction $A$, i.e., $\mu(A) = \underset{t \to 0^{+}}{\lim} \frac{\|I+tA\|-1}{t}$.

See, for example, Table 1 of \cite{maidens2014reachability} for commonly used definition of the induced matrix measure. For example, with the infinity norm, $\mu_{\infty}(A):=\max_i \left(a_{ii} + \sum_{j\neq i} |a_{ij}|\right)$. 

Observe that $\mu(A)$ can be negative. In fact, one of the properties of the matrix measure is, for all eigenvalues $\lambda_i(A)$ of $A$, 
\begin{equation} \label{eq:mu_mat}
    \mathbf{R}(\lambda_i(A)) \leq \mu(A) \leq \|A\|,
\end{equation}
where $\mathbf{R}(\lambda)$ extracts the real part of a complex number $\lambda$. 
 
\section{Problem formulation} \label{sec:problem}
Consider a continuous-time nonlinear system 
\begin{equation} \label{eq:sys}
    \dot{x}=f(x,d), \qquad x(0)\in K,
\end{equation}
where the state is $x\in\mathbb{R}^{n_x}$, the input $d\in\mathbb{R}^{n_d}$ is seen as a disturbance which resides in a closed set $D$ that contains the origin. This input signal $d$ can be any measurable, locally essentially bounded function of time to the set $D$, and the set of all such inputs is denoted as $\mathcal{D}$. The function $f:\mathbb{R}^{n_x}\times \mathbb{R}^{n_d} \to\mathbb{R}^{n_x}$ is $C^1$ (continuously differentiable) with respect to the first argument and $f(0,0)=0$. The set $K$ is compact and convex. We denote the solution to \eqref{eq:sys} initialised at $x(0)=x_0$ evaluated at time $t$ with input $d$ as $x(t,x_0,d)$. We assume that system \eqref{eq:sys} is forward complete, i.e., the solutions to $x(t,x_0,d)$ exist for all time $t\geq 0$.

We are interested in system \eqref{eq:sys} with solutions that converge to each other, other than converging to an equilibrium point, at a desired convergence rate $\beta\in\mathcal{KL}$ for a particular time horizon $[0,T]$, $T>0$ and when the states are initialized from a compact set $K$. 

\begin{defn} \label{def:pres_stability}
For a given $\beta\in\mathcal{KL}$ and initial set $K$, system \eqref{eq:sys} is 
\begin{itemize}
    \item \textit{$(\beta, K)$-incrementally asymptotically stable}, if for all $d\in\mathcal{D}$, for any $x_1$, $x_2\in K$, the solutions to system \eqref{eq:sys} satisfy the following for all $t\geq 0$,  $$|x(t,x_1,d)-x(t,x_2,d)| < \beta(|x_1-x_2|,t).$$
    
    \item  \textit{$(\beta, K, \varepsilon)$-incrementally practically stable}, for  $\varepsilon>0$, if for any $d_1$, $d_2\in\mathcal{D}$, any $x_1$, $x_2 \in K$, the solutions to system \eqref{eq:sys} satisfy the following for all $t\geq 0$,  $$|x(t,x_1, d_1)-x(t,x_2,d_2)|  < \beta(|x_1-x_2|,t) + \varepsilon.$$ 
    \item When $\beta(r,t)=e^{-\alpha t} r$, for $\alpha>0$, we say that system \eqref{eq:sys} is $(\alpha,K)$-\textit{incrementally exponentially stable} or $(\alpha,K,\varepsilon)$\textit{-incrementally practically exponentially stable}, respectively. \hfill $\Box$
\end{itemize}
\end{defn}

The prescribed stability notions in Definition \ref{def:pres_stability} are incremental stability properties, cf. \cite{zamani2011backstepping}, where these notions are highly relevant for state estimation. 

The \textit{date rate} (bit rate) of system \eqref{eq:sys} corresponds to the number of samples of the state space constrained to the initial set $K$ that are needed per unit time to achieve the prescribed stability notions in Definition \ref{def:pres_stability}. By letting the time tend to infinity, the notion of \textit{entropy} introduced in the Section \ref{sec:entropy} captures the growth rate of these numbers.  This paper provides bounds on the entropy of a system \eqref{eq:sys} under disturbances that possesses incremental stability properties in Section \ref{sec:bounds}. The computation of these bounds is achieved through the construction of approximating sets in Section \ref{sec:finite_sets}, which induced a robust state estimation algorithm for systems under disturbances using quantized and time-sampled state measurements in Section \ref{sec:algo}.

\section{Estimation entropy} \label{sec:entropy}
We call the finite set of points $S=\{x_1,x_2,\dots,x_N\} \subset K$, $(T,\beta,K,\varepsilon)$-asymptotically approximating or $(T,\beta,K,\varepsilon)$-practically approximating when the corresponding notions of stability in Definition \ref{def:pres_stability} are achieved within a finite time interval $[0,T]$. We formalise this below.

\begin{defn} \label{def:approx_def}
Let $T>0$, $\beta\in\mathcal{KL}$, $\varepsilon>0$ and $m\in\mathbb{N}_{\geq 1}$. Given a finite set of points $S=(x_1,x_2,\dots,x_n)\subset K$, system \eqref{eq:sys} is 
\begin{itemize}
    \item $(T,\beta,K,m\varepsilon)$-\textit{asymptotically approximating}, if for any initial state $x_0\in K$, any $d \in \mathcal{D}$, there exists a point $x_i\in S$ such that the following holds for all $ t\in[0,T]$,
\begin{equation} \label{eq:asym_approx}
    |x(t,x_0,d)-x(t,x_i,d)|< \beta(|x_0-x_i|+m\varepsilon,t).
\end{equation} 
    \item $(T,\beta,K,m\varepsilon)$-\textit{practically approximating}, if for any initial state $x_0\in K$, any $d_0$, $d_i \in \mathcal{D}$, there exists a point $x_i\in S$ such that the following holds for all $ t\in[0,T]$,
\begin{align} \label{eq:prac_approx}
    |x(t,&x_0,d_0)-x(t,x_i,d_i)| \nonumber \\ &< \beta(|x_0-x_i|+m\varepsilon,t)+m\varepsilon.
\end{align} 
\item When the class $\mathcal{KL}$ function $\beta(r,t):=e^{-\alpha t}r$, where $\alpha>0$, we say that the sets are $(T,\alpha,K,m\varepsilon)$-\textit{exponentially} or \textit{practically exponentially approximating}, respectively.\hfill $\Box$
\end{itemize}

\end{defn}

\begin{rem} \label{rem:def_exist}
    In Definition \ref{def:approx_def}, $\varepsilon>0$ is introduced in \eqref{eq:asym_approx} to ensure the existence of a finite asymptotically approximating set, which we will see later in Section \ref{sec:finite_sets}. This was also observed in the case for practical stabilization in \cite{colonius2021entropy}.
\end{rem}

Let $s_a(T,\beta,K,\varepsilon)$ and $s_p(T,\beta,K,\varepsilon)$ denote the minimal cardinality of the $(T,\beta,K,\varepsilon)$-asymptotically approximating set and the $(T,\beta,K,\varepsilon)$-practically approximating set, respectively. Suppose that the approximating sets are finite, we define the the asymptotic and practical estimation entropy respectively, below.

\begin{defn} \label{def:entropy}
Let $T>0$, $\beta\in\mathcal{KL}$ and $\varepsilon>0$.
\begin{itemize}
    \item The asymptotic estimation entropy is 
    \begin{equation} \label{eq:ae_entropy}
        h_a(\beta,K) := \underset{\varepsilon\to 0}{\lim} \underset{T\to\infty}{\limsup} \, \frac{1}{T} \log s_a(T,\beta,K,\varepsilon).
    \end{equation}
    \item The $\varepsilon$-practical estimation entropy is
        \begin{equation} \label{eq:pe_entropy}
        h_{\varepsilon p}(\beta,K,\varepsilon) := %\underset{\varepsilon\to 0}{\lim}
        \underset{T\to\infty}{\limsup} \, \frac{1}{T} \log s_p(T,\beta,K,\varepsilon).
    \end{equation}
%    \item The practical estimation entropy is 
%    \begin{equation} \label{eq:p_entropy}
%        h_{p}(\beta,K) :=\underset{\varepsilon\to 0}{\lim} \, h_{\varepsilon p}(\beta,K,\varepsilon).
%    \end{equation}
    \item When $\beta(r,t)=e^{-\alpha t}r$, where $\alpha >0$, the estimation entropies $h_a(\beta,K)$ and $h_{\varepsilon p}(\alpha,K,\varepsilon)$ are denoted as $h_a(\beta,K)$ and $h_{\varepsilon p}(\alpha,K,\varepsilon)$, respectively. \hfill $\Box$
\end{itemize} 
\end{defn}

Since $s_a$ and $s_p$ are the minimum number of functions needed to approximate the state $x$ with the desired accuracy and convergence rate, the notions of estimation entropy $h_a$ and $h_p$ are the respective average number of state estimates over the time interval $[0,T]$. The growth rate of $s_a$ or $s_p$ over time is captured by the inner $\limsup$ and by taking its limit as $\varepsilon$ goes to 0, we obtain the worst case over $\varepsilon>0$. Equivalently, the outer limit coincides with the supremum over $\varepsilon$. We use the natural logarithm here for convenience since we consider encoded information taking values in $\mathbb{R}^{n_x}$ (the same choice was made in \cite{colonius2021entropy}). This is in place of the usual choice of the logarithm in base $2$ due to its direct relation to the number of 0's and 1's needed in the encoded information that is transmitted over a digital communication channel as considered in \cite{liberzon2017entropy}, \cite{matveev2019observation}, \cite{katok1997introduction}, \cite{downarowicz2011entropy}, for example.  

\section{Existence of finite approximating sets} \label{sec:finite_sets}
We follow the Bowen \cite{bowen1971entropy} and Dinaburg \cite{dinaburg1971relations} construction on entropy via spanning sets, which are called approximating sets (Definition \ref{def:approx_def}) in this paper. In this section, we investigate the existence of \textit{finite} approximating sets for both incremental asymptotic and practical stability properties. This forms an important intermediary step towards obtaining entropy bounds later.

The following notion of a \textit{$\delta$-cover} is useful in quantifying the resolution of the bit rate, which we recall here. For a bounded set $K\subseteq \mathbb{R}^{n}$ and $\delta > 0$, a $\delta$-cover is a finite collection of points $S=\{x_1,x_2,\dots,x_N\}$ such that the set $K$ is a subset of the union of closed balls centered a $x_i$ with radius $\delta$, i.e., $K\subseteq \bigcup_{i\in\mathbb{N}_{[1,N]}} B(x_i,\delta)$.

Central to approximating the $\delta$-cover of $K$ for the approximating sets is a result from \cite{maidens2014reachability} that allows us to bound the distance between its trajectories as a function of their initial distance and its convergence rate, is given by the matrix measure of the system's Jacobian in $x$, which is uniformly bounded as stated below.
\begin{assum} \label{assum:mat_bound}
Consider system \eqref{eq:sys}, where there exists $\bar{\mu}\in\mathbb{R}$ such that
    $\mu\left(\frac{\partial f}{\partial x}(x,d)\right) \leq \bar{\mu}$, for all $(x,d)\in K \times \mathcal{D}.$ 
\end{assum}
\begin{prop}[Prop. 1 in \cite{maidens2014reachability}] \label{prop:mat_close}
    Consider system \eqref{eq:sys} under Assumption \ref{assum:mat_bound}. Then, for any initial state $x_1$, $x_2\in K$ and any $d\in\mathcal{D}$, the corresponding solutions of \eqref{eq:sys} satisfy
    \begin{equation}
        |x(t,x_1,d)-x(t,x_2,d)| \leq e^{\bar{\mu}t} |x_1-x_2|, \; t\in[0,T].
    \end{equation} \hfill $\Box$
\end{prop}

 When $\bar{\mu} < 0$, then system \eqref{eq:sys} is incrementally exponentially stable \cite{angeli2002lyapunov}. On the other hand, if $f$ of system \eqref{eq:sys} is $C^1$, a positive $\bar{\mu}$ provides a sharper bound\footnote{Examples where the rate was exploited instead of the Lipschitz constant of $f$ are the model detection algorithm of \cite{liberzon2017entropy} and in the approximation of reachability sets in \cite{maidens2014reachability}.} over the Lipschitz constant $L$ of $f$ due to \eqref{eq:mu_mat}, where the induced norm of the Jacobian of $f$ in $x$ is equal to the Lipschitz constant $L$.    

We are now ready to establish the existence of finite approximating sets and its resolution for $(\beta,K)$-incrementally stable and $(\beta,K,\varepsilon)$-incrementally practically stable systems, as well as their exponential counterparts (the special case where $\beta(r,t):=e^{-\alpha t}r$ for $\alpha >0$), respectively. The results are detailed in Lemmas \ref{lem:finite_practical}-\ref{lem:finite_stable_exp} below and summarised in Table \ref{tab:finite_summary}.

\begin{lem}[incremental practical stability] \label{lem:finite_practical}
Consider system \eqref{eq:sys} under Assumption \ref{assum:mat_bound}. Given  $\beta\in\mathcal{KL}$ and $\varepsilon>0$, suppose for any $x_0\in K$, any $d_1, d_2 \in\mathcal{D}$, there exists $x\in K$ such that for all $t\geq 0$,
\begin{equation} \label{eq:prac_x_cond}
    |x(t,x_0,d_1)-x(t,x,d_2)|< \beta(|x_0-x|,t)+\varepsilon.
\end{equation}
Then, for all $T>0$, there exists a finite collection of points $S:=\{x_1,x_2,\dots,x_N\}\subset K$ such that for any $x_0\in K$, any $d_1$, $d_2\in \mathcal{D}$, there exists a point $x_i \in S$ such that for $t\in[0,T]$,
\begin{equation} \label{eq:prac_x_conc}
    |x(t,x_0,d_1)-x(t,x_i,d_2)|<\beta(|x_0-x_i|+\varepsilon,t) + 2\varepsilon.
\end{equation}
Moreover, 
$S$ is a $\delta$-cover of $K$, where $\delta:=e^{-\bar{c}T} \varepsilon$ with  $\bar{c}:=\max\{\bar{\mu},0\}$.
\hfill $\Box$
\end{lem}
\begin{lem}[incremental practical exponential stability] \label{lem:finite_practical_exp}
Consider system \eqref{eq:sys} under Assumption \ref{assum:mat_bound}. Given  $\alpha>0$ and $\varepsilon>0$, suppose for any $x_0\in K$, any $d_1, d_2 \in\mathcal{D}$, there exists $x\in K$ such that for all $t\geq 0$,
\begin{equation} \label{eq:exp_prac_x_cond}
    |x(t,x_0,d_1)-x(t,x,d_2)|< e^{-\alpha t}|x_0-x|+\varepsilon.
\end{equation}
Then, for all $T>0$, there exists a finite collection of points $S:=\{x_1,x_2,\dots,x_N\}$ such that for any $x_0\in K$, any $d_1$, $d_2\in \mathcal{D}$, there exists a point $x_i \in S$ such that for $t\in[0,T]$,
\begin{equation} \label{eq:exp_prac_x_conc}
    |x(t,x_0,d_1)-x(t,x_i,d_2)|<e^{-\alpha t}(|x_0-x_i|+\varepsilon) + \varepsilon.
\end{equation}
Moreover, 
$S$ is a $\delta$-cover of $K$, where $\delta:=e^{-(M_{p}+\alpha)T} \varepsilon$ with  $M_p:=\max \{\bar{\mu},-\alpha\}$.
\hfill $\Box$
\end{lem}

\begin{lem}[incremental stability] \label{lem:finite_stable}
Consider system \eqref{eq:sys} under Assumption \ref{assum:mat_bound}. Given  $\beta\in\mathcal{KL}$, suppose for any $x_0\in K$, any $d \in\mathcal{D}$, there exists $x\in K$ such that for all $t\geq 0$,
\begin{equation} \label{eq:stable_x_cond}
    |x(t,x_0,d)-x(t,x,d)|< \beta(|x_0-x|,t).
\end{equation}
Then, for all $T>0$ and $\varepsilon>0$,  there exists a finite collection of points $S:=\{x_1,x_2,\dots,x_N\}\subset K$ such that for any $x_0\in K$, any  $d\in \mathcal{D}$, there exists a point $x_i \in S$ such that for $t\in[0,T]$,
\begin{equation} \label{eq:stable_x_conc}
    |x(t,x_0,d)-x(t,x_i,d)|<\beta(|x_0-x_i|+\varepsilon,t) + \varepsilon.
\end{equation}
Moreover, 
$S$ is a $\delta$-cover of $K$, where $\delta:=e^{-\bar{c}T} \varepsilon$ with  $\bar{c}:=\max\{\bar{\mu},0\}$.
\hfill $\Box$
\end{lem}
\begin{lem}[incremental exponential stability] \label{lem:finite_stable_exp}
Consider system \eqref{eq:sys} under Assumption \ref{assum:mat_bound}. Given  $\alpha>0$ and $\varepsilon>0$, suppose for any $x_0\in K$, any $d \in\mathcal{D}$, there exists $x\in K$ such that for all $t\geq 0$,
\begin{equation} \label{eq:exp_stable_x_cond}
    |x(t,x_0,d)-x(t,x,d)|< e^{-\alpha t}|x_0-x|.
\end{equation}
Then, for all $T>0$ and $\varepsilon>0$,  there exists a finite collection of points $S:=\{x_1,x_2,\dots,x_N\}\subset K$ such that for any $x_0\in K$, any  $d\in \mathcal{D}$, there exists a point $x_i \in S$ such that for $t\in[0,T]$,
\begin{equation} \label{eq:exp_stable_x_conc}
    |x(t,x_0,d)-x(t,x_i,d)|<e^{-\alpha t}(|x_0-x_i|+2\varepsilon).
\end{equation}
Moreover, 
$S$ is a $\delta$-cover of $K$, where $\delta:=e^{-(M_{e}+\alpha)T} \varepsilon$ with  $M_{e}:=\max\{\bar{\mu},-\alpha\}$.
\hfill $\Box$
\end{lem}

\begin{table}[h!]
    \centering
    \begin{tabular}{|p{0.08\columnwidth}||p{0.41\columnwidth}|p{0.29\columnwidth}|}
    \hline
    Lemma & Results for a given $\beta\in\mathcal{KL}$, $\alpha>0$ compact $K\subset \mathbb{R}^{n}$, $\varepsilon>0$ and $T>0$. & Resolution $\delta$ \\ \hline \hline
     \ref{lem:finite_practical} & $(\beta,K,\varepsilon)$-prac. incrementally stable  $\implies$ \par finite $(T,\beta,K,2\varepsilon)$-prac. approximating set & \par $e^{-\bar{c}T}\varepsilon$, \par $\bar{c}:=\max\{\bar{\mu},0\}.$ \\ \hline
     \ref{lem:finite_practical_exp} & $(\alpha,K,\varepsilon)$-prac. incrementally exp. stable  $\implies$ \par finite $(T,\alpha,K,\varepsilon)$-prac. exp. approximating set & \par $e^{-(M+\alpha)T}\varepsilon$, \par $M:=\max\{\bar{\mu},-\alpha\}.$ \\ \hline
     \ref{lem:finite_stable} & $(\beta,K)$- incrementally stable  $\implies$  \par finite $(T,\beta,K,\varepsilon)$-prac. approximating set &  \par $e^{-\bar{c}T}\varepsilon$,  \par $\bar{c}:=\max\{\bar{\mu},0\}.$ \\ \hline
     \ref{lem:finite_stable_exp} & $(\alpha,K)$-incrementally exp. stable  $\implies$ \par finite $(T,\alpha,K,2\varepsilon)$-exp. approximating set & \par $e^{-(M+\alpha)T}\varepsilon$, \par $M:=\max\{\bar{\mu},-\alpha\}.$ \\ \hline
\end{tabular}
    \caption{Summary of Lemmas \ref{lem:finite_practical}-\ref{lem:finite_stable_exp} for finite approximating sets} \vspace{-2em}
    \label{tab:finite_summary}
\end{table}

From Lemmas \ref{lem:finite_practical}-\ref{lem:finite_stable_exp} (cf. Table \ref{tab:finite_summary}), we see that a consequence of sampling the state space within a finite time interval is a decrease in approximation accuracy. 

\begin{rem}
We can also obtain a finite  approximating set $S$ in the absence of Assumption \ref{assum:mat_bound}, if $f$ is locally Lipschitz in $x$ on $[0,T]\times \mathcal{D}$, with Lipschitz constant $L$. In this case, we establish closeness of solutions for system \eqref{eq:sys} using the Gronwall-Bellman inequality (cf. \cite[Theorem 3.4]{khalil2002nonlinear}). By the same way the proof of Lemma \ref{lem:finite_practical} was done, we obtain that $S$ is a $e^{-LT} \varepsilon$-cover of the compact set $K$. However, under Assumption \ref{assum:mat_bound} and by \eqref{eq:mu_mat} where the induced matrix norm of the Jacobian of $f(x,d)$ in $x$ is the Lipschitz constant of $f$, we see that Assumption \ref{assum:mat_bound} provides a sharper bound. \hfill $\Box$
\end{rem}

%We quantify the growth rate of the average number of state samples or quantized states which are needed via computing bounds on the estimation entropy in the next section.

\section{Estimation entropy bounds} \label{sec:bounds}
We obtain the upper and lower bounds on the various types of estimation entropy (see Definition \ref{def:entropy}) in this section.

\begin{thm}[Upper bounds] \label{thm:upper}
Consider system \eqref{eq:sys} under Assumption \ref{assum:mat_bound}. Given $\beta\in\mathcal{KL}$, $\varepsilon>0$ and $\alpha>0$, suppose system \eqref{eq:sys} is 
\begin{enumerate}[(i)]
    \item $(\beta,K,\varepsilon)$-practically incrementally stable. Then, the $2\varepsilon$-practical estimation entropy satisfies $h_{\varepsilon p}(\beta,K,2\varepsilon) \leq \bar{c} n$, where  $\bar{c}:=\max\{\bar{\mu},0\}$ comes from Lemma \ref{lem:finite_practical}.
    \item $(\alpha,K,\varepsilon)$-practically incrementally exponentially stable. Then, the $\varepsilon$-practical estimation entropy satisfies $h_{\varepsilon p}(\alpha,K,\varepsilon) \leq (M_p+\alpha) n$, where $M:=\max\{\bar{\mu},-\alpha\}$ comes from Lemma \ref{lem:finite_practical_exp}.
    \item $(\beta,K)$-incrementally stable. Then, the $\varepsilon$-practical estimation entropy satisfies
    $h_{\varepsilon p}(\alpha,K,\varepsilon) \leq \bar{c} n$, where $\bar{c}:=\max\{\bar{\mu},0\}$ comes from Lemma \ref{lem:finite_stable}.
    \item $(\alpha,K)$-incrementally exponentially stable. Then, the asymptotic estimation entropy satisfies $h_a(\alpha,K)\leq (M+\alpha)n$, where $M:=\max\{\bar{\mu}, -\alpha\}$ comes from Lemma \ref{lem:finite_stable_exp}. \hfill $\Box$
\end{enumerate}    
    
\end{thm}

Next, we derive a lower bound on the estimation entropy. Following a volume growth argument employed in both \cite[Proposition 3]{liberzon2017entropy} for estimation entropy of systems without inputs and \cite[Theorem 3.4]{colonius2021entropy} for stabilization entropy, we compute the lower bound on the estimation entropy for systems with inputs.

In the result that follows, we denote the set of all solutions of system \eqref{eq:sys} in the time interval $[0,T]$, $T<\infty$ with initial states in $K$ and disturbance $d\in\mathcal{D}$  by $R_{K,\mathcal{D}}[0,T]:=\{x(t,x_0,d) \, | \, x_0 \in K, d\in \mathcal{D}, t\in[0,T] \}$.

\begin{thm}[Lower bounds] \label{thm:lower}
Consider system \eqref{eq:sys} under Assumption \ref{assum:mat_bound}.   
\begin{enumerate}[(i)]
    \item Given $\beta\in\mathcal{KL}$ and $\varepsilon>0$, the $\varepsilon$-practical estimation entropy satisfies \begin{align}\label{eq:lower_ent_ep} 
        h_{\varepsilon p}(\beta,K,\varepsilon)\geq \underset{(x,d)\in R_{K,\mathcal{D}}[0,T]\times \mathcal{D}}{\inf}\textrm{tr} f_{x}(x, d ). 
    \end{align}
    \item Given $\beta(r,t):=e^{-\alpha t} r$ and $\varepsilon>0$, where $\alpha>0$, the $\varepsilon$-practical estimation entropy satisfies
    \begin{align} \label{eq:lower_ent_ep_exp}
        h_{\varepsilon p}(\alpha,K,\varepsilon)  \geq & \, (M_p+\alpha)n \nonumber \\
        &+ \underset{(x,d)\in R_{K,\mathcal{D}}[0,T]\times \mathcal{D}}{\inf}\textrm{tr} f_{x}(x, d ), 
    \end{align} 
    and the asymptotic estimation entropy satisfies
    \begin{align} \label{eq:lower_ent_a_exp}
        h_{a}(\alpha,K)\geq &\, (M_e+\alpha)n \nonumber \\
        &+ \underset{(x,d)\in R_{K,\mathcal{D}}[0,T]\times \mathcal{D}}{\inf}\textrm{tr} f_{x}(x, d ),
    \end{align}
    where $M_p=M_e:=\max\{\bar{\mu},-\alpha\}$. \hfill $\Box$
\end{enumerate} 
\end{thm}

In the case where $\textrm{tr} f_x(x,d)$ is unbounded from below over $R_{K,\mathcal{D}}[0,T]$, the lower bound of the estimation entropies are $-\infty$. Further, if $\textrm{tr} f_x(x,d)$ is negative, the estimation entropies can be negative, which is not a useful lower bound. We expect that sharper bounds can be derived with the method in \cite{kawan2021remote}.

\section{Robust state estimation } \label{sec:algo}
\begin{figure}[h!]
    \centering
    \includegraphics[width=9cm]{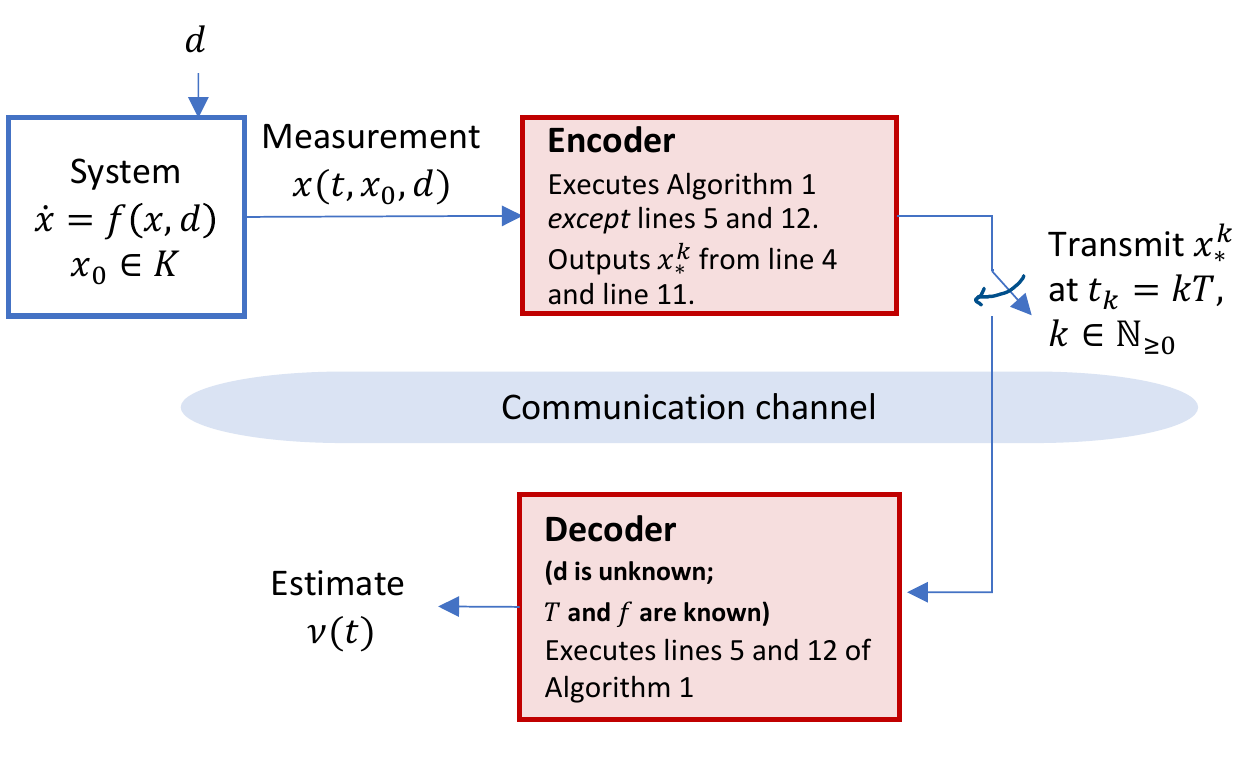} \vspace{-2em}
    \caption{State estimation with quantized and time-sampled measurements}
    \label{fig:setup}
\end{figure}

The construction of the finite approximation sets in Section \ref{sec:finite_sets} induces a state estimation algorithm for a system under disturbance \eqref{eq:sys} using quantized and time-sampled measurements. At discrete instances in time $t_k=kT$, $T\in(0,\infty)$, $k\in\mathbb{N}_{\geq 0}$, a limited amount of data about system \eqref{eq:sys} is transmitted. The exact details of the transmitted data is depicted in Figure \ref{fig:setup} and explained below. The decoder receives the transmitted data and produces an estimate $\nu(t)$ of the state $x$ of system \eqref{eq:sys} under the conditions where the time interval $T$ between transmissions is known, the dynamics $f$ of system \eqref{eq:sys} are known, but not the disturbance $d$. This setup is depicted in Figure \ref{fig:setup}. This encoding-decoding setup is standard practice for control and estimation over a finite capacity communication channel, see \cite{nair2007feedback} for an overview. \vspace{-1em}

The robust state estimation algorithm produces an estimate $\nu:\mathbb{R}_{\geq 0}\to \mathbb{R}^{n_x}$ which converges to a neighborhood of the true state $x$ of system \eqref{eq:sys}, where the neighbourhood is dependant on the accuracy of the approximator $\Sigma_a$ that is used on the decoder side. We assume that system \eqref{eq:sys} satisfies Assumption \ref{assum:mat_bound} and is $(\alpha,K,\varepsilon)$-practically incrementally stable (Definition \ref{def:pres_stability}) for constructive purposes. Hence, the approximator $\Sigma_a$ is chosen to be $\dot{\hat{x}}=f(\hat{x},0)$ with an initial condition $x_*^k$ specified by Algorithm \ref{alg:cap}. In fact, this initial point is the data that is transmitted over the communication channel.

From here on, all norms are infinity norms and $B(r,\delta)$ balls are hypercubes centered at $r$ with diameter $2\delta$.

The robust state estimation algorithm proposed in this section resembles the algorithm presented in \cite{liberzon2017entropy} with a modification (see line 8 in Algorithm \ref{alg:cap})  to accommodate for estimation inaccuracies due to the presence of disturbances. We describe the procedure in Algorithm \ref{alg:cap}.

\begin{algorithm}
\caption{Robust state estimation algorithm with quantized and time-sampled measurements}\label{alg:cap}
\begin{algorithmic}[1]
\Require $\alpha$, $\varepsilon$, $K$, $T$, $M$, approximator $\Sigma_a$, sequence of measurements $\{x(kT,x_0,d)\}_{k\in\mathbb{N}_{\geq 0}}$..  

\State Initialise: $c_0$ is the center of $K$, $d_0=\textrm{diameter}(K)/2$, $\delta_0=d_0e^{-(M+\alpha)T}$.

\State Construct $K_0=B(c_0,d_0)$. 
\State Generate $\delta_0$-cover of $K_0$, denoted by $S_0$ which has $N_0$ points.
\State Given measurement $x(0,x_0,d)$, choose $x^0_*\in S_0$ which satisfies
$|x(0,x_0,d)-x_*^0|\leq d_0 e^{-(M+\alpha)T}.$
\State Compute $\nu_{0}(t)=x(t,x_*^0,0)$, for $t\in[0,T]$ using approximator $\Sigma_a$ initialized at $x_*^0$.
\State $k=1$.
\For {$k\in\mathbb{N}_{\geq 1}$}
\State Update $\delta_{k}=e^{-\alpha T} \delta_{k-1}+2\varepsilon$.
\State Construct $K_k=B(\nu_{k-1}(T),\delta_k)$.
\State Generate $\delta_k e^{-(\alpha+M)T}$-cover of $K_k$, denoted by $S_k$ which has $N_k$ points.
\State Given measurement $x(kT,x_0,d)$, choose $x^k_*\in S_k$ which satisfies
$|x(kT,x_0,d)-x_*^k|\leq \delta_k e^{-(M+\alpha)T}$.
\State Compute $\nu_{k}(t)=x(t,x_*^k,0)$, for $t\in[0,T]$ using approximator $\Sigma_a$ initialized at $x_*^k$. 
\State Increment $k$ by 1.
\EndFor
\end{algorithmic}
\end{algorithm}

At initialization $t=0$, the algorithm receives the following information:
\begin{itemize}
    \item initial set $K$,
%    \item number of quantization samples $N$,
    \item duration in between each transmission $T$,
    \item vector field $f$ of the system \eqref{eq:sys}, but not disturbance $d$,
    \item the constant $M$ from Lemma \ref{lem:finite_practical_exp}, and
    \item an approximator $\Sigma_a$ to compute the solutions $x(t,z_0,0)$ for any $z_0\in K$, on the decoder side.
\end{itemize}

The initial parameters are $K_0 = B(c_0,d_0)$, where $c_0$ and $d_0$ are the center and diameter of the hypercube which are chosen such that $K\subseteq K_0$.   A $d_0 e^{-(M+\alpha) T}$-cover of $K_0$ is generated,  which we denote by $S_0:=\{x^{0}_{1},x^{0}_{2},\dots,x^{0}_{N_0}\}$. Here, $M$ comes from Lemma \ref{lem:finite_practical_exp}. Given that the measurement $x(0,x_0,d)$ is available at $t=0$, we choose ball $B(x^{0}_{*},d_0 e^{-(M+\alpha) T})$, where  $x^{0}_{*}=x^{0}_{i} \in S_0$ which contains the measurement $x(0,x_0,d)$\footnote{Although we have the true initial condition $x_0$ of system \eqref{eq:sys}, we cannot reconstruct its state $x(t,x_0,d)$ on the decoder side as the disturbance $d$ is unknown.}. Using the point $x_*^0$ which is closest to the current measurement $x(0,x_0,d)$, we generate an approximating function $\nu_0(t)=x(t,x_*^0,0)$ for $t\in[0,T]$.

At each iteration $k\in\mathbb{N}_{\geq 1}$, we compute the following:
\begin{enumerate}
    \item the geometrically shrinking radius $\delta_k=e^{-\alpha T}\delta_{k-1} + 2\varepsilon$ of $K_k$, where the presence of $2\varepsilon$ is to account for the approximation inaccuracy $\varepsilon>0$ which we know from Lemma \ref{lem:finite_practical_exp},
    \item the hypercube $K_k$ which is centered at the last approximate $\nu_{k-1}(T)$ with radius $\delta_k$,
    \item the $\delta_k e^{-(M+\alpha)T}$-cover of $K_k$, denoted by $S_k:=\{x_1^{k},x_2^{k},\dots,x_{N_k}^{k} \}$,
    \item the closest point $x_*^{k}=x_i^{k}\in S_k$ to the received measurement $x(kT,x_0,d)$ and generate the approximating function $\nu_{k}(t)=x(t,x_*^{k},0)$ for $t\in[0,T]$.
\end{enumerate}

Under the operating conditions we have described, 
we show that the current measurement $x(kT,x_0,d)$ is always contained within the current set $K_k$ and hence, the closest point $x^k_*$ that is used to generate the next set $K_{k+1}$ is always within the ball $B(x(kT,x_0,d),\delta_k e^{-(M+\alpha)T})$. This property will be used to show the accuracy of the state estimate $\nu(t)$, which we define to be a piecewise constant signal that concatenates all the approximating functions $v_k(t)$ for $t\in[0,T]$, as follows
\begin{equation} \label{eq:nu}
    \nu(t) = \nu_{k}(t-kT), \; t\in[kT,(k+1)T).
\end{equation}

We are now ready to provide the following guarantees about the robust state estimation algorithm.

\begin{thm} \label{thm:algo}
Consider system \eqref{eq:sys} under Assumption \ref{assum:mat_bound} and is $(\alpha,K,\varepsilon)$-incrementally practically exponentially stable. Then, for all $T>0$, Algorithm \ref{alg:cap} satisfies the following for all $k\in\mathbb{N}_{\geq 0}$
\begin{enumerate}[(i)]
    \item $x(kT,x_0,d)\in K_k$, and
    \item $|x(t,x_0,d)-\nu(t)|\leq e^{-\alpha t} d_0 + 2\varepsilon (\bar{a}_{k}+1),$ for all $t\in[kT,(k+1)T)$,  where $\bar{a}_{k}:=\frac{1-a^k}{1-a}$, $a=e^{-\alpha T}$, and
    \item $\underset{t\to\infty}{\limsup}\;|x(t,x_0,d)-\nu(t)| \leq 2\varepsilon\left(({1-e^{-\alpha T}})^{-1}+1 \right).$ 
\end{enumerate} \hfill $\Box$
\end{thm}
When system \eqref{eq:sys} under Assumption \ref{assum:mat_bound} is $(\alpha,K)$-incrementally exponentially stable, the proof of Theorem \ref{thm:algo} holds true for $\varepsilon=0$, resulting in exponential convergence of the state estimate.  
\begin{cor}
Consider system \eqref{eq:sys} under Assumption \ref{assum:mat_bound} and is $(\alpha,K)$-incrementally exponentially stable. Then, for all $T>0$, Algorithm \ref{alg:cap} with $\varepsilon=0$ guarantees (i) of Theorem \ref{thm:algo} and $|x(t,x_0,d)-\nu(t)|\leq e^{-\alpha t}d_0$ for all $t\in[kT,(k+1)T)$, $k\in\mathbb{N}_{\geq 1}$. \hfill $\Box$
\end{cor}
Due to space constraints, the data rate of Algorithm \ref{alg:cap} in the setup depicted in Figure \ref{fig:setup} and its relation to entropy are not discussed here, and will be investigated in the future.

\section{Conclusions and future work} \label{sec:conclude}
We have introduced notions of estimation entropy for $\varepsilon$-incremental practical stability and incremental stability of deterministic continuous-time systems under disturbances. Bounds on the estimation entropies are obtained via approximating sets. Finally, a robust state estimation algorithm is proposed for systems under disturbances using quantized and time-sampled measurements. Future work include investigating the asymptotic estimation entropy; and relating the $\varepsilon$-practical estimation entropy to minimal data rates of a finite capacity communication channel, to name a few. 
%%%%%%%%%%%%%%%%%%%%%%%%%%%%%%%%%%%%%%%%%%%%%%%%%%%%%%%%%%%%%%%%%%%%%%%%%%%%%%%%

%%%%%%%%%%%%%%%%%%%%%%%%%%%%%%%%%%%%%%%%%%%%%%%%%%%%%%%%%%%%%%%%%%%%%%%%%%%%%%%%
\appendix
\subsection{Proof of Lemma \ref{lem:finite_practical}}
For any $x_0\in K$ and any $d_1,d_2\in\mathcal{D}$, choose an $x\in K$ such that for all $t\in[0,T]$,
\begin{align} \label{eq:choose_x}
    |x(t,x_0,d_1)-x(t,x,d_2)|&< \beta(|x_0-x|,t)+\varepsilon.
\end{align}
Using the triangle inequality, we obtain
\begin{align} \label{eq:xd_1}
    &|x(t,x_0,d_1) -x(t,x_i,d_2)| \nonumber \\
    & \leq |x(t,x_0,d_1)-x(t,x,d_2)| + |x(t,x,d_2)-x(t,x_i,d_2)| \nonumber \\
    & < \beta(|x_0-x|,t)+\varepsilon + |x(t,x,d_2)-x(t,x_i,d_2)| \nonumber \\
    & \leq \beta(|x_0-x_i|+|x_i-x|,t)+\varepsilon \nonumber \\ 
    & \qquad+ |x(t,x,d_2)-x(t,x_i,d_2)|.
\end{align}
Under Assumption \ref{assum:mat_bound}, we apply Proposition \ref{prop:mat_close} to obtain that 
\begin{align} \label{eq:xd_2}
    |x(t,x,d_2)-x(t,x_i,d_2)| \leq e^{\bar{\mu}t}|x-x_i|, \; t\in[0,T].
\end{align}

We proceed by constructing the set $S=\{x_1,x_2,\dots,x_N\} \subset K$ and by showing that it is a $\delta$-cover of the set $K$ ($\delta>0$ will be constructed below), as well as a $(T,\beta,K,2\varepsilon)$-practically approximating set (i.e., for any $x_0\in K$, \eqref{eq:prac_x_conc} holds for at least one $x_i\in S$). 

To this end, the points $x_i\in S$ are points $x_i \in K$ chosen such that it is a $\delta$-cover of $K$, where $\delta:=e^{-\bar{c}T} \varepsilon$, with $\bar{c}=\max\{\bar{\mu},0\}\geq 0$ as given. In other words, we have constructed the set $S=\{x_1,x_2,\dots,x_N\} \subset K$.

By this construction, we see that for the chosen $x\in K$, we can always choose an $x_i \in S$ such that $|x-x_i|< e^{-\bar{c}T} \varepsilon$. Then, from \eqref{eq:xd_1} and \eqref{eq:xd_2}, we obtain
\begin{align}
    &|x(t,x_0,d_1) -x(t,x_i,d_2)| \nonumber \\
    & < \beta(|x_0-x_i|+e^{-\bar{c}T} \varepsilon,t)+\varepsilon + e^{\bar{\mu}t}e^{-\bar{c}T} \varepsilon
\end{align}
Since $\bar{c}\geq 0$, $e^{-\bar{c}T} \leq e^{-\bar{c}t}$ for $t\leq T$, and $\bar{c}-\bar{\mu}\geq 0$. Thus $e^{\bar{\mu}t} e^{-\bar{c}T} \leq e^{\bar{\mu}t} e^{-\bar{c}t} \leq  e^{-(\bar{c}-\bar{\mu})t} \leq 1$. Therefore, we obtain \eqref{eq:prac_x_conc}.\hfill $\Box$

\subsection{Proof of Lemma \ref{lem:finite_practical_exp}}
Given $\alpha>0$, the proof takes similar steps as the proof of Lemma \ref{lem:finite_practical} from \eqref{eq:choose_x} to \eqref{eq:xd_2}, by taking $\beta(r,t)=e^{-\alpha t} r$. 

The point of departure is in the construction of the set $S:=(x_1,x_2,\dots,x_N)\subset K$ to be a $\delta$-cover of $K$ and to show that it is a $(T,\alpha,K,\varepsilon)$-exponentially approximating set. To this end, we choose $x_i \in S$ such that $S$ is a $\delta$-cover of $K$, where $\delta:=e^{-(M+\alpha)T}\varepsilon$. Then, for any $x\in K$, we can always choose a point $x_i\in S$ such that $|x-x_i|<e^{-(M+\alpha)T}\varepsilon$. 

From \eqref{eq:xd_1} and \eqref{eq:xd_2}, we have for all $t\in[0,T]$,
\begin{align} \label{eq:exp_p_1}
    |x(t,x_0,d_1) -x(t,x_i,d_2)| < & e^{-\alpha t}  |x_0-x_i| \nonumber \\
    &+ e^{-(\alpha-\bar{\mu}) t} |x-x_i|. 
\end{align}

Since $M+\alpha \geq 0$, we have that $e^{-(M+\alpha)T} \leq e^{-(M+\alpha)t}$ for all $t\in[0,T]$. Hence,  $e^{-(\alpha-\bar{\mu}) t} |x-x_i| \leq  e^{-(\alpha-\bar{\mu}+M+\alpha) t} \leq e^{-2\alpha t} \varepsilon$. From \eqref{eq:exp_p_1}, we thus obtain \eqref{eq:exp_prac_x_conc}. \hfill $\Box$

\subsection{Proofs of Lemma \ref{lem:finite_stable} and \ref{lem:finite_stable_exp}}
The proof of Lemma \ref{lem:finite_stable}  (Lemma \ref{lem:finite_stable_exp}) proceed in a similar manner as the proof for Lemma \ref{lem:finite_practical} (Lemma \ref{lem:finite_practical_exp}).

\subsection{Proof of Theorem \ref{thm:upper}}
We will perform the proof for case (i). The other cases can be proven in the same manner by employing Lemma \ref{lem:finite_practical_exp} for case (ii), Lemma \ref{lem:finite_stable} for case (iii) and Lemma  \ref{lem:finite_stable_exp} for case (iv), respectively.

\underline{Proof for case (i)}:
Since system \eqref{eq:sys} is $(\beta,K,\varepsilon)$-incrementally practically stable, we apply Lemma \ref{lem:finite_practical} to obtain a finite $(T,\beta,K,2\varepsilon)$-practically approximating set $S:=\{x_1,x_2,\dots,x_N\}$ which is also a $\delta$-cover of the compact set $K$, where $\delta:=e^{-\bar{c}T}\varepsilon$, with  $\bar{c}>|\bar{\mu}|$. The rest of the proof follows the same mechanism as in the proof for Proposition 2 in \cite{liberzon2017entropy}.

Let $c(\delta,K)$ denote the minimum cardinality of a bounded set $K$ with balls of radius $\delta$. Then, we have that $s_p(T,\beta,K,\varepsilon) \leq c(e^{-\bar{c}T}\varepsilon,K)$, where we recall that $s_p(T,\beta,K,2\varepsilon)$ is the minimum cardinality of the $S$. Then, a bound on the $2\varepsilon$-practical estimation entropy $h_{\varepsilon p}(T,\beta,K,2\varepsilon)$ is
\begin{align}
    \underset{T\to\infty}{\limsup} \, \frac{1}{T} \log s_p(T,\beta,K,2\varepsilon) & \leq \underset{T\to\infty}{\limsup} \, \frac{1}{T} \log c(e^{-\bar{c}T}\varepsilon,K) \nonumber \\
    & \leq \underset{T\to\infty}{\limsup} \, \bar{c} \log  \frac{c(e^{-\bar{c}T}\varepsilon,K)}{\bar{c}T}. \nonumber
\end{align}
By noting that $\bar{c}T=\log(e^{\bar{c}T}\varepsilon^{-1}\cdot \varepsilon) = \log(e^{\bar{c}T}\varepsilon^{-1}) + \log \varepsilon$, we obtain
\begin{align} \label{eq:upper_1}
    & \underset{T\to\infty}{\limsup} \, \frac{1}{T} \log s_p(T,\beta,K,2\varepsilon) \nonumber \\
    & \leq \underset{T\to\infty}{\limsup} \, \bar{c}   \frac{\log c(e^{-\bar{c}T}\varepsilon,K)}{ \log(e^{\bar{c}T}\varepsilon^{-1}) + \log \varepsilon } \nonumber \\
    & \leq \bar{c} \, \underset{T\to\infty}{\limsup} \,   \frac{\log c(e^{-\bar{c}T}\varepsilon,K)}{ \log(e^{\bar{c}T}\varepsilon^{-1})},
\end{align}
where the last bound was obtained as $\log \, \varepsilon$ does not affect the limit superior. Recall the fact that the upper box dimension (c.f. \cite[Chapter 2.2]{falconer2004fractal}) satisfies the following for any bounded set $K \subset \mathbb{R}^{n}$, $    \underset{T\to\infty}{\limsup} \, \frac{\log \, c(\delta,K)}{\log(\delta^{-1})} \leq n.$
Therefore, we conclude that $h_{\varepsilon p}(\beta,K,2\varepsilon) \leq \bar{c} n$. 

\subsection{Proof of Theorem \ref{thm:lower}}
\underline{We first prove (i)}. For a $(T,\beta,K,\varepsilon)$-practically approximating set $S=\{x_1,x_2,\dots,x_N\}$, recall from Lemma \ref{lem:finite_stable} that the $R_{K,\mathcal{D}}[0,T]$ is covered by balls  $B(x(t,x_i,d),\delta(T,\varepsilon))$ for $x_i\in S$, $d\in\mathcal{D}$, $t\in[0,T]$, where $\delta(T,\varepsilon):=\beta(e^{-\bar{c}T}\varepsilon+\varepsilon,0)+\varepsilon$ is obtained from \eqref{eq:prac_x_conc} of Lemma \ref{lem:finite_practical}. Hence, the minimal cardinality $s_p(T,\beta,K,\varepsilon)$ of $S$ is lower bounded by the volume of $R_{K,\mathcal{D}}[0,T]$ and the volume of each ball $B(x(t,x_i,d),\delta(T,\varepsilon))$, i.e.,
\begin{align} \label{eq:vol_gen}
    s_p(T,\beta,K,\varepsilon) \geq \frac{\textrm{vol}\left(R_{K,\mathcal{D}}[0,T]\right)}{\textrm{vol}\left( B(x(t,x_i,d),\delta(T,\varepsilon)) \right)}.
\end{align}
We first derive a lower bound for $\textrm{vol}\left(R_{K,\mathcal{D}}[0,T]\right)$. Let $z=R_{K,\mathcal{D}}[0,T]$ and observe that
\begin{align} \label{eq:vol_nume}
\textrm{vol}&\left(R_{K,\mathcal{D}}[0,T]\right):=\int_{R_{K,\mathcal{D}}[0,T]} dz = \int_{K} |\textrm{det}\, x_{x}(T,x,d)| dx, \nonumber \\
& \geq \underset{(x,d)\in K\times \mathcal{D}}{\inf} \textrm{vol}(K) \cdot |\textrm{det}\, x_{x}(T,x,d)| 
\end{align}
where $x_x = \frac{\partial x(T,x,d)}{\partial x}$ and the second equality is obtained by a change of integration variables. By the Abel-Jacobi-Liouville identity \cite[Lemma 3.11]{teschl2012ordinary}, we have
\begin{equation} \label{eq:liouville}
    \textrm{det}\, x_{x}(t,x,d) = \textrm{exp}\left(\int_{0}^{t} \textrm{tr} f_{x}(x(s,x,d), d(s) )\,ds \right)
\end{equation}
where $f_x(x,d)=\frac{\partial f(x,d)}{\partial x}$. Hence, from \eqref{eq:liouville} and \eqref{eq:vol_nume}, we obtain the following
\begin{align}
    &\textrm{vol}\left(R_{K,\mathcal{D}}[0,T]\right) \nonumber \\
    &\geq \underset{(x,d)\in K\times \mathcal{D}}{\inf} \textrm{vol}(K) \cdot \textrm{exp}\left(\int_{0}^{t} \textrm{tr} f_{x}(x(s,x,d), d(s) )\,ds \right) \nonumber \\
    & =  \textrm{vol}(K) \cdot \textrm{exp}\left(\underset{(x,d)\in K\times \mathcal{D}}{\inf}\int_{0}^{t} \textrm{tr} f_{x}(x(s,x,d), d(s) )\,ds \right) \nonumber \\
    & \geq \textrm{vol}(K) \cdot \textrm{exp}\left(T \underset{(x,d)\in R_{K,\mathcal{D}}[0,T]\times \mathcal{D}}{\inf}\textrm{tr} f_{x}(x, d ) \right)
\end{align}

Using the fact that $\textrm{vol}\left( B(x(t,x_i,d),\delta(T,\varepsilon)) \right)=(2\delta(T,\varepsilon))^n$ by employing the infinity norm and from \eqref{eq:vol_gen} and \eqref{eq:vol_nume}, we obtain  \begin{equation}\label{eq:sp_almost}    s_p(T,\beta,K,\varepsilon) \geq \frac{\textrm{vol}(K) \cdot \textrm{exp}\left(T \underset{(x,d)\in K\times \mathcal{D}}{\inf}\textrm{tr} f_{x}(x, d ) \right)}{(2\delta(T,\varepsilon))^n}.
\end{equation} Hence,
\begin{align}
    \frac{1}{T} \log s_p(T,\beta,K,\varepsilon) \geq &\frac{1}{T}  \log \left( {\textrm{vol}(K)} \right) - \frac{1}{T} \log  {(2\delta(T,\varepsilon))^n}  \nonumber \\
    & + \left(\underset{(x,d)\in R_{K,\mathcal{D}}[0,T]\times \mathcal{D}}{\inf}\textrm{tr} f_{x}(x, d )\right). \nonumber
\end{align}
By taking the $\limsup_{T\to\infty}$ for the inequality above, we observe that $\limsup_{T\to\infty}\frac{1}{T}  \log \left( \frac{\textrm{vol}(K)}{(2\delta(T,\varepsilon))^n} \right)=0$ as $\delta(T,\varepsilon)\to \beta(\varepsilon,0)+\varepsilon$ as $T\to\infty$, where we recall that $\delta(T,\varepsilon):=\beta(e^{-\bar{c}T}\varepsilon+\varepsilon,0)+\varepsilon$. Therefore, we obtain \eqref{eq:lower_ent_ep}.

\underline{Case (ii)} can be proven in the same manner as case (i). For case (ii), we apply Lemma \ref{lem:finite_practical_exp} and obtain a $(T,\alpha,K,\varepsilon)$-practically exponentially approximating set $S_{pe}:=\{x_1,x_2,\dots,x_n\}$ such that $R_{K,\mathcal{D}}[0,T]$ is covered by balls $B(x(t,x_i,d),\delta_{pe}(T,\varepsilon))$, where $\delta_{pe}(T,\varepsilon):=e^{-(M_p+\alpha)T}\varepsilon+\varepsilon$, which we obtain from \eqref{eq:exp_prac_x_conc} of Lemma \ref{lem:finite_practical_exp}. 

The point of departure from the proof for case (i) lies in the radius of the balls covering $R_{K,\mathcal{D}}[0,T]$. Hence, we obtain
\begin{align}
    \underset{T\to\infty}{\limsup} &\frac{1}{T} \log(2 \delta_{pe}(T,\varepsilon))^{n} & \nonumber \\
    & = \underset{T\to\infty}{\limsup} \frac{1}{T} \log(2 e^{-(M_p+\alpha)T}\varepsilon+2\varepsilon)^{n} \nonumber \\
    & = \underset{T\to\infty}{\limsup} \frac{1}{T} \log(2 e^{-(M_p+\alpha)T}\varepsilon)^{n} \nonumber \\
    & = \underset{T\to\infty}{\limsup} \frac{1}{T} \log( e^{-(M_p+\alpha)T})^{n} + \underset{T\to\infty}{\limsup} \frac{1}{T} \log(2 \varepsilon )^{n} \nonumber \\
    & = -(M_p+\alpha)n, \label{eq:sp_del}
\end{align}
where we obtain the second equality as the additive $\varepsilon$ term does not affect the $\limsup$ as $T\to\infty$. Finally, from \eqref{eq:sp_almost} and \eqref{eq:sp_del}, we obtain \eqref{eq:lower_ent_ep_exp} as desired.

\underline{Case (iii)} can be proven \textit{mutatis mutandis} by application of Lemma \ref{lem:finite_stable_exp}. \hfill $\Box$

\subsection{Proof of Theorem \ref{thm:algo}}
\underline{Proof of (i)}: At $k=0$, by the construction of $K_0$ [line 2 of Algorithm \ref{alg:cap}], the initial set $K\subseteq K_0$. Since $x(0,x_0,d) =x_0 \in K$, we have that $x(0,x_0,d)\in K_0$.
We now show that $x(kT,x_0,d)\in K_k$ for $k\in\mathbb{N}_{\geq 1}$, or that $x(kT,x_0,d)\in B(\nu_{k-1}(T),\delta_k)$ [line 9 of Algorithm \ref{alg:cap}].

First, note that $x(kT,x_0,d)=x(T,x((k-1)T,x_0,d),d)$ and $\nu_{k-1}(T)=x(T,x_*^{k-1},0)$. Hence,
\begin{align}
    |x(kT,&x_0,d) - \nu_{k-1}(T)| \nonumber \\ &=|x(T,x((k-1)T,x_0,d),d)-x(T,x_*^{k-1},0)|. \nonumber 
\end{align}
Since system \eqref{eq:sys} under Assumption \ref{assum:mat_bound} is $(\alpha,K,\varepsilon)$-incrementally practically exponentially stable, we apply Lemma \ref{lem:finite_practical_exp} and obtain
\begin{align}
    |x&(kT,x_0,d) - \nu_{k-1}(T)| \nonumber \\ & < e^{-\alpha T}\left(|x((k-1)T,x_0,d)-x_*^{k-1}|+\varepsilon \right)+\varepsilon \nonumber \\
    & \leq e^{-\alpha T}\left(\delta_{k-1}e^{-(M+\alpha)T}+\varepsilon \right)+\varepsilon \nonumber \\
    & < e^{-\alpha T}\left(\delta_{k-1}+\varepsilon \right)+\varepsilon < e^{-\alpha T}\delta_{k-1} + 2 \varepsilon=\delta_k, \label{eq:proof_algo}
\end{align}
where we obtained the second inequality from line 11 of Algorithm \ref{alg:cap}, and the third and fourth inequality due to $e^{-(M+\alpha)T}\leq 1$ and $e^{-\alpha T} < 1$ since $\alpha >0$ and $M+\alpha \geq 0$. From \eqref{eq:proof_algo}, we have shown that $x(kT,x_0,d)\in B(\nu_{k-1}(T),\delta_k)$ and thus (i) is satisfied.

\underline{Proof of (ii)}: Observe that for $t\in[kT,(k+1)T]$, $x(t,x_0,d)=x(t-kT,x(kT,x_0,d),d)$ and $\nu(t)=\nu_{k}(t-kT)=x(t-kT,x_*^k,0)$. Therefore,
\begin{align}
    |x&(t,x_0,d)-\nu(t)| \nonumber \\
    &=|x(t-kT,x(kT,x_0,d),d)-x(t-kT,x_*^k,0)| \nonumber \\
    & < e^{-\alpha(t-kT)}(|x(kT,x_0,d)-x_*^k|+\varepsilon) + \varepsilon \nonumber   
\end{align}
where we obtained the last inequality thanks to Lemma \ref{lem:finite_practical_exp} since system \eqref{eq:sys} under Assumption \ref{assum:mat_bound} is $(\alpha,K,\varepsilon)-$incrementally practically exponentially stable.
From line 11 of Algorithm \ref{alg:cap}, we have that $|x(kT,x_0,d)-x_*^k|\leq \delta_k e^{-(M+\alpha)T}$ and by solving line 8 of Algorithm \ref{alg:cap} iteratively, we obtain $\delta_{k}=a^{k}\delta_0 + 2\varepsilon \bar{a}_{k}$, where $\bar{a}_{k}:=\frac{1-a^{k}}{1-a}$ with $a:=e^{-\alpha T}$ as defined in the theorem. Therefore,
\begin{align}
    |x&(t,x_0,d)-\nu(t)| \nonumber \\
    & < e^{-\alpha t} e^{\alpha k T}(e^{-\alpha k T}\delta_0+2\varepsilon \bar{a}_{k}+\varepsilon) + \varepsilon \nonumber \\
    & \leq e^{-\alpha t} \delta_0 + 2 \varepsilon \bar{a}_{k} + 2 \varepsilon,
\end{align}
since $e^{-\alpha(t-kT)}\leq 1$ for $t\in[kT,(k+1)T]$. Thus, we have proven (ii).

\underline{Proof of (iii)}: Since $\underset{t\to\infty}{\limsup} \; e^{-\alpha t} \delta_0 = 0$ and $\underset{t\to\infty}{\limsup}\; 2 \varepsilon \bar{a}_{k} = \underset{k\to\infty}{\limsup} \; 2 \varepsilon \bar{a}_{k} = 2\varepsilon \frac{1}{1-a}$ as $a<1$, we have shown (iii). \hfill $\Box$

\bibliographystyle{ieeetr}
\bibliography{global_bib.bib}

\end{document}